\documentclass[12pt]{article}

\usepackage{latexsym}
\usepackage{amsmath}
\usepackage{amssymb}
\usepackage{amsfonts}

\def\disp{\displaystyle}
\def\o{\over}

\def\reals{I\!\!R}
\def\a{\alpha}
\def\ve{\varepsilon}
\def\g{\gamma}
\def\gm{\gamma(dy,du)}
\def\G{\Gamma}
\def\U{{\cal U}}
\def\P{{\cal P}}

\def\ph{\varphi}
\def\l{C}
\def\d{\delta}

\def\bu{\bar u}
\def\by{\bar y}
\def\ve{\varepsilon}

\def\hf{\hfill{$\Box$}}

\def\bY{\bar Y}
\def\bG{\bar G}

\def\bW{\bar W}

\def\l{\lambda}

\def\gm{\,d\g}

\newtheorem{Theorem}{Theorem}[section]
\newtheorem{Proposition}[Theorem]{Proposition}

\newtheorem{Lemma}[Theorem]{Lemma}

\begin{document}

\begin{center}\Large{{\bf Compactification Method in Linear Programming Approach to Infinite-Horizon Optimal Control Problems with a Noncompact State Constraint}}
\end{center}

{\bf Ilya Shvartsman}\\
{\it \small{Department of Computer Science and Mathematics,
Penn State Harrisburg, Middletown, PA 17057, USA}}

\bigskip


\bigskip

{\bf Abstract.} This paper is devoted to a study of infinite horizon optimal control problems with time discounting and time averaging criteria in discrete time. It is known that these problems are related to certain infinite-dimensional linear programming problems, but compactness of the state constraint is a common assumption imposed in analysis of these LP problems.  In this paper, we consider an unbounded state constraint and use Alexandroff compactification to carry out the analysis. We also establish asymptotic relationships between the optimal values of problems with time discounting and long-run average criteria.

\section{Introduction}\label{Intro}

In this paper, we consider discrete time controlled dynamical system
\begin{equation}\label{A1}
\begin{aligned}
&y(t+1)=f(y(t),u(t)), \; t=0,1,\dots\, \\
&y(0)=y_0,\\
&y(t)\in Y,\\
&u(t)\in U(y(t)).
\end{aligned}
\end{equation}
Here  $Y$ is a given nonempty closed subset of $\reals^m$, $\ U(\cdot):\,Y\leadsto U_0$ is an upper semicontinuous compact-valued mapping to a given compact metric space $U_0$,
$\ f(\cdot, \cdot ):\,\reals^m\times U_0\to \reals^m$ is a continuous function.  

A control $u(\cdot)$ and the pair $(y(\cdot),u(\cdot))$ are called an {\em admissible control} and an {\em admissible process}, respectively, if the relationships \eqref{A1} are satisfied. The sets of admissible controls are denoted by $\U(y_0)$ or $\U_S(y_0)$, depending on whether the problem is considered on the  infinite time horizon or on a finite time sequence  ($t\in \{0,\dots,S-1\}$, where $S$ is a positive integer).

Throughout the text we assume that $\U(y_0)\neq \emptyset$ for all $y_0\in Y$. (Systems that satisfy this property are called {\em viable} on $Y$.)

On the trajectories of \eqref{A1} we consider the optimal control problems
\begin{equation}\label{A111}
(1-\a)\inf_{u(\cdot)\in \U(y_0)}\sum_{t=0}^{\infty} \a^t g(y(t),u(t))=:V_{\a}(y_0),
\end{equation} 
and
\begin{equation}\label{A112-1}
{1\o S}\min_{u(\cdot)\in \U_S(y_0)}\sum_{t=0}^{S-1} g(y(t),u(t))=:V(S, y_0)
\end{equation}
where $g:\,\reals^m\times U_0\to \reals^m$ is a continuous function bounded below   and  $\a\in (0,1)$ is a discount factor. 

The limits $\disp \lim_{\a\to 1^-} V_{\a}(y_0)$ and $\disp \lim_{S\to \infty} V(S,y_0)$ are called {\em Abel and Ces\`aro limits}, respectively. Limits of {Abel} and  Ces\`aro types have been studied in various contexts. There are a lot of publications dedicated to the conditions ensuring their existence and equality in problems of dynamic programming and optimal control in discrete and continuous time, see, e.g., \cite{BQR-2015,GS-Opt-22,Sorin92,OV-2012,ShvPAFA21}. It is shown in \cite{OV-2012} in continuous time setting that the limits $\disp \lim_{\a\to 1^-} V_{\a}(y_0)$ and $\disp \lim_{S\to \infty} V(S,y_0)$ must be equal when the convergence is uniform. This condition is replaced in \cite{GS-Opt-22} with a weaker condition that the limits are continuous functions of $y_0$.

The {\em linear programming approach} to problems of control of nonlinear dynamical systems is based on the fact that
the so-called {\em occupational measures} generated by admissible processes satisfy certain 
{\em linear} equations that represent the dynamics of the system in integral form. This makes it
possible to reformulate optimal control problems as infinite-dimensional linear
programming (IDLP) problems considered on the spaces of occupational measures. Solutions of
the dual problems to these IDLP problems can be used to construct feedback controls that ensure
optimality of the corresponding trajectories. This approach has been used in many works in deterministic and stochastic settings in continuous and discrete time, see, e.g., \cite{BG}, \cite{GQ}-\cite{GS-DCDS-22},\cite{ShvCOT}.  For other approaches for dealing with discrete time optimal control problems on infinite time horizon see the survey \cite{HLerma23} and references therein. 

A significant limitation of the linear programming approach is the usual assumption that the state constraint $Y$ is compact, which is needed to ensure compactness of the space of probability measures $\P(Y\times U_0)$ in weak$^*$ topology. Compactness of $Y$ in not assumed in \cite{BG05}, but in \cite{BG05} an assumption on the occupational measures is made that is difficult to verify. In this paper, we obtain generalizations of some of the results from  \cite{GPS-2017} {\em without the compactness assumption} on $Y$. The significance of these results is elaborated on at the end of the next section. 

The paper is organized as follows. In Section \ref{Prelim}, we formulate problems \eqref{A111} and \eqref{A112-1} in terms of the occupational measures and state known results used in the sequel. In Section \ref{Compact} we use Alexandroff compactification to compactify $Y$ and formulate the corresponding problems in the new spaces. The main results of the paper that include asymptotic representation of the sets of occupational measures and the limits $\disp \lim_{\a\to 1^-} V_{\a}(y_0)$ and $\disp \lim_{S\to \infty} V(S,y_0)$ are in Sections \ref{S-a} and \ref{S-b}.

\section{Occupational Measures and Preliminaries}\label{Prelim}

Denote
$$
G:=\{(y,u)|\,y\in Y,\,u\in U(y),\,f(y,u)\in Y\}.
$$
Any admissible process stays in $G$ for all $t$.

In the paper \cite{GPS-2017}, problems \eqref{A111} and \eqref{A112-1} are considered under the assumption that the set $Y$ is {\em compact}.
In \cite{GPS-2017}, the space of probability measures on Borel subsets of $G$ is denoted by $\P(G)$ and the sets below are defined:
\begin{equation}\label{e-W1}
\begin{aligned}
&W_{\a}(y_0):=\{\g\in \P(G)|\, \\
&\int_{G}[\a(\ph(f(y,u))-\ph(y))
+(1-\a)(\ph(y_0)-\ph(y))]\,\gm=0\quad \forall \ph\in C(Y)\},
\end{aligned}
\end{equation}
\begin{equation}\label{M17}
\begin{aligned}
W:=\{\g\in \P(G)|\, \int_{G}(\ph(f(y,u))-\ph(y))\,
\gm=0\quad \forall \ph\in C(Y)\}.
\end{aligned}
\end{equation}
The definition of $W$ formally follows from the definition of $W_{\a}(y_0)$ by setting $\a=1$. Consider the following 
infinite-dimensional linear programming (IDLP) problems:
\begin{equation}\label{D1}
\min_{\g\in W_{\a}(y_0)} \int_{G} g(y,u)\,d\g=:g^*_{\a}(y_0),
\end{equation}
\begin{equation}\label{M22}
\min_{\g\in W} \int_{G} g(y,u)\,d\g=:g^*.
\end{equation}
(These problems are problems of linear programming because the objective functions and the constraints are linear in the ``decision variable" $\g$.)
It is proved in \cite{GPS-2017} (see  formulas (8) and (9) in \cite{GPS-2017}) that 
\begin{equation}\label{eq-Res1}
V_{\a}(y_0)=g^*_{\a}(y_0)
\end{equation}
and
\begin{equation}\label{eq-Res2}
\lim_{\a\to 1^-}\min_{y\in Y}V_{\a}(y)=\lim_{S\to \infty}\min_{y\in Y}V(S,y)=g^*.
\end{equation}
(Note that in \cite{GPS-2017} $V_{\a}(y_0)$ is defined without the factor $(1-\a)$ in \eqref{A111} and $V(S,y_0)$ is defined without the factor $1/S$.)

Problems \eqref{A111} and \eqref{A112-1} can be equivalently formulated in terms of the so-called {\em occupational measures}. For an admissible process $(y(\cdot),u(\cdot))$, a probability measure $\g^{\a}_{u}$ is called the {\em
 discounted occupational measure generated by $u(\cdot)$} if, for any Borel set $Q\subset G$,
\begin{equation}\label{EE6}
\g^{\a}_{u}(Q)=(1-\a) \sum_{t=0}^{\infty} \a^t 1_Q(y(t),u(t)),
\end{equation}
where $1_Q(\cdot)$ is the indicator function of $Q$. 
One can see that this definition is equivalent to the relationship
\begin{equation}\label{G8}
\int_{G} q(y,u) \g^{\a}_u(dy,du)=(1-\a)\sum_{t=0}^{\infty} \a^t q(y(t),u(t))
\end{equation}
for any Borel measurable function $q$ on $G$. Indeed, \eqref{EE6} obviously implies \eqref{G8} for a function which is a finite sum of indicator functions of Borel measurable sets, and  the validity of \eqref{G8} for an arbitrary Borel $q$  follows from the definition of the Lebesgue integral as a limit of integrals of simple functions; see, e.g. \cite{Ash}. 

Similarly, a probability measure
$\g_{u,S}$ is called the {\em averaging occupational measure generated by $u(\cdot)$} over the time sequence $\{0,1,...,S-1 \}$ if, for any Borel set $Q\subset G$,
\begin{equation}\label{G81}
\g_{u,S}(Q)={1\o S} \sum_{t=0}^{S-1} 1_Q(y(t),u(t)).
\end{equation}
Let $\G_{\a}(y_0)$ and $\G(S,y_0)$ denote the sets of all discounted and averaging occupational measures generated by the admissible controls, that is,
$$
\G_{\a}(y_0):=\bigcup_{u(\cdot)\in \U(y_0)}\g^{\a}_{u} \;\hbox{ and }\,  \G(S,y_0):=\bigcup_{u(\cdot)\in \U_S(y_0)}\g_{u,S}.
$$
Note that due to the assumption that $Y$ is viable, the sets $\G_{\a}(y_0)$ and $\G(S,y_0)$ are not empty for all $y_0$ and $S$. Due to \eqref{G8} and \eqref{G81},  problems \eqref{A111} and \eqref{A112-1}
can be rewritten as
\begin{equation}\label{CV11}
\inf_{\g\in \G_{\a}(y_0)} \int_{G} g(y,u)d\g= V_{\a}(y_0)
\end{equation}
and 
\begin{equation*}
\inf_{\g\in \G(S,y_0)} \int_{G} g(y,u)d\g= V(S,y_0).
\end{equation*}

Denote
$$
\G_{\a}:=\bigcup_{y_0\in Y}(\G_{\a}(y_0)) \,\hbox{ and }\, \G(S):=\bigcup_{y_0\in Y}(\G(S,y_0)).
$$
It proved in Theorem 5.4 in \cite{GPS-2017} that in the case if $Y$ is compact, the limits of the closure of convex hulls of $\G_{\a}$ and $\G(S)$ exist and 
\begin{equation}\label{CV1}
\lim_{\a\to 1^-}\bar{\rm co}\,\G_{\a}=\lim_{S\to \infty}\bar{\rm co}\,\G(S)=W,
\end{equation}
where convergence is understood in weak$^*$ sense. 

The proofs of \eqref{eq-Res2} and \eqref{CV1} in \cite{GPS-2017} essentially rely on boundedness of $Y$, which implies compactness of the set  $W$ in weak$^*$ topology. The goal of this paper is to obtain generalizations of \eqref{eq-Res2} and \eqref{CV1} to the case when the set $Y$ is closed and {\em unbounded}. 

Equality \eqref{CV1} plays an important role in obtaining estimates of the limits of $V_{\a}(y_0)$ as $\a\to 1^-$ and of $V(S,y_0)$ as $S\to \infty$ in \cite{BGS} in a non-ergodic case, that is, when these limits depend on the initial condition. Relation \eqref{CV1} is also used, in a continuous-time framework, for deriving a representation formula of the cluster points in the uniform convergence topology of the functions $V_{\a}(\cdot)$ and $V(S,\cdot)$ (see \cite{BQR-2015}). Furthermore, \eqref{CV1} is used  in the analysis of singularly perturbed control systems (see \cite{G04}). The significance of the results obtained in this paper is that they open the door for carrying over the corresponding results to the case of unbounded state space, which will be a subject of further research.

\section{Compactification}\label{Compact}

Throughout the rest of the paper we assume that $Y$ is an unbounded closed set.

Recall that the one-point compactification (or Alexandroff compactification) of a locally compact noncompact topological space $X$ is $\bar X=X\cup\{\infty\}$, where open sets on $\bar X$ consist of open sets of $X$ and sets of the form $(X\setminus C)\cup\{\infty\}$, where $C$ is closed and compact in $X$.

Let $\bY:=Y\cup\{\infty_Y\}$ and $\bG:=G\cup\{\infty_G\}$ be one-point compactifications of $Y$ and $G$, respectively. Since $U_0$ is compact, one can see that $(y,u)\to \infty_G$ if and only if $y\to \infty_Y$. For this reason, we can identify $\infty_G$ with $(\infty_Y,u)$ for any $u\in U_0$.

Denote by $C(\bY)$ the space of continuous bounded functions  on $\bY$. It follows from the topology of $\bY$ that function $\ph:\bY\to \reals$ belongs to $C(\bY)$ if $\ph$  is continuous on $Y$ and has a finite limit as $|y|\to \infty$. In this case, we set $\ph(\infty_Y):=\disp \lim_{|y|\to \infty} \ph(y).$

Similarly, $q:\,\bG\to \reals$ belongs to $C(\bG)$ if $q$ is continuous on $G$, and  the finite limit $\disp q(\infty_G):=\lim_{(y,u)\to \infty_G} q(y,u)$ exists. The latter is equivalent to the existence of $\disp\lim_{|y|\to \infty} q(y,u)$, its independence of $u$ and uniformity with respect to $u$.  

To extend the state space of \eqref{A1} to include the ``infinite" state, set 
\begin{equation}\label{f-inf}
f(\infty_Y,u)=\infty_Y \;\forall u\in U_0.
\end{equation}
This implies that  if $y_0=\infty_Y$ then $y(t)=\infty_Y$ for all $t$. Conversely, the only trajectory that reaches $\infty_Y$ is the one that starts at $\infty_Y$. 

To be able to extend $f$ to $\infty_G$ by continuity, that is, to ensure that $\disp \lim_{(y,u)\to \infty_G}f(y,u)=\infty_Y$, throughout the paper we assume that
\begin{equation}\label{f-inf2}
|f(y,u)|\to \infty \hbox{ as }|y|\to \infty \quad\hbox{uniformly in }u.
\end{equation}
This assumption is not too restrictive. It holds, for example, for the linear system $y(t+1)=Ay(t)+Bu(t)$ if $A$ is not singular;  for a difference system $y(t+1)=y(t)+v(y(t),u(t))$, it is implied by the condition $|v(y,u)|\le c(1+|y|^{a})$ for some $c$ and $0<a<1$, for all $u$.

To preserve upper semicontinuity of the map $U(\cdot)$, set $U(\infty_Y):=U_0$. By $\P(\bG)$ we denote the space of probability measures on Borel subsets of $\bG$.

Due to the Riesz-Markov-Kakutani theorem, the space dual to $C(\bG)$ is the space of regular Borel measures on $\bG$. Due to the Banach-Alaoglu theorem, the unit ball in $(C(\bG))^*$ is weakly$^*$ compact. These theorems imply the following. 
\begin{Proposition}\label{P-comp}
The space  $\P(\bG)$ is weakly$^*$ compact, that is, if $\g_k\in  \P(\bG)$ is a sequence of probability measures, then there exist a subsequence $\g_{k'}$ and $\g\in  \P(\bG)$ such that for any $q(y,u)\in C(\bG)$
$$
\lim_{k'\to \infty} \int_{\bG} q(y,u)\,d\g_{k'}=\int_{\bG} q(y,u)\,d\g.
$$
\end{Proposition}

{\bf Proof.} The existence of $\g\in(C(\bG))^*$ such that of $\g_k$ converges to $\g$ weakly$^*$ along a subsequence  follows from weak$^*$ compactness of the unit ball in $(C(\bG))^*$. To show that $\g\in \P(\bG)$ it only remains to prove that $\g(\bG)=1$. Indeed,
$$
\g(\bG)=\int_{\bG}1\,d\g=\lim_{k'\to \infty}\int_{\bG}1\,d\g_{k'}=1.
$$
\hf

\section{Generalization of \eqref{CV1}}\label{S-a}

Let
\begin{equation}\label{CV31}
W:=\{\g\in \P(G)|\, \int_{G}(\ph(f(y,u))-\ph(y))
\gm=0\quad \forall \ph\in C(\bY)\}.
\end{equation}
(Note that in the case of unbounded $Y$ we take $\ph\in C(\bY)$, as opposed to $\ph\in C(Y)$ in \eqref{M17} when $Y$ is bounded.)
We also define a subset $\bW$ of $\P(\bG)$ as a ``compact counterpart of $W$": 
\begin{equation}\label{CV4}
\bW:=\{\g\in \P(\bG)|\, \int_{\bG}(\ph(f(y,u))-\ph(y))
\gm=0\quad \forall \ph\in C(\bY)\}.
\end{equation}
For the integral in \eqref{CV4} to be well defined, the integrand must have a finite limit as $(y,u)\to \infty_G$. Due to \eqref{f-inf2}, for any $\ph\in C(\bY)$ we have $\disp\lim_{(y,u)\to \infty_G}(\ph(f(y,u))-\ph(y))=0$, so this condition is satisfied.

Since $\bW$ is a subset of $\P(\bG)$, it may contain measures $\g$ such that $\g(\{\infty_G\})>0$;
$W$ doesn't contain such measures being a subset of $\P(G)$. Since $\bW$ is a closed subset of a weakly$^*$ compact set $\P(\bG)$, $\bW$ is weakly$^*$ compact. 

Below we write $\g_k\to \g$ if $\g_k$ converges to $\g$ weakly$^*$ on $\P(\bG)$.

From \eqref{CV31} and \eqref{CV4} it follows that 
\begin{equation}\label{W-W}
W=\bW\cap \P(G).
\end{equation}
The following proposition gives a condition under which $W$ is not empty. (Hence, $\bW$ is also not empty due to \eqref{W-W}.)

\begin{Proposition}\label{P-nonempty}
If there exists an admissible process $(y(\cdot),u(\cdot))$ that remains in a bounded subset of $G$ for all $t$, then $W\neq \emptyset$.
\end{Proposition}
{\bf Proof.} Let $(y(\cdot),u(\cdot))\subset D$ for all $t$, where $D$ is a closed bounded subset of $G$. Take a sequence $\a_k\to 1^-$. For each $k$, this process generates a discounted occupational measure $\g_{k}:=\g_u^{\a_k}$ supported on $D$. Along a subsequence (we do not relabel) we have $\g_k\to \g\in \P(\bG)$ . Let us show that $\g$ is supported on $D$. Indeed, let $q(y,u)=\min\{\hbox{dist}\, ((y,u),D),1\}$. This function belongs to $C(\bG)$, is equal to zero on $D$ and is positive outside of $D$. We have
$$
0=\int_{D}q\,d\g_k=\int_{\bG}q\,d\g_k\to \int_{\bG}q\,d\g.
$$
If $\g$ had support outside of $D$, then the last integral would be positive. Therefore, $\g$ is supported on $D$.

Denote 
$$
W_{\a}(y_0):=\{\g\in \P(G)|\, 
\int_{G}[\a(\ph(f(y,u))-\ph(y))
+(1-\a)(\ph(y_0)-\ph(y))]\gm=0\quad \forall \ph\in C(\bY)\}.
$$
Since $\g_k\in \G_{\a_k}(y_0)$ and $\G_{\a_k}(y_0)\subset W_{\a_k}(y_0)$, where $y_0$ is the initial state of $(y(\cdot),u(\cdot))$ (see \cite{GPS-2017}, Proposition 5 for the proof of the latter inclusion), we have 
\begin{equation*}
\begin{aligned}
0=&\int_{G}[\a_k(\ph(f(y,u))-\ph(y))+(1-\a_k)(\ph(y_0)-\ph(y))]\,d\g_k\\
=&\int_{\bG}[\a_k(\ph(f(y,u))-\ph(y))+(1-\a_k)(\ph(y_0)-\ph(y))]\,d\g_k\\
=&\int_{\bG}[(\a_k-1)(\ph(f(y,u))-\ph(y))+(1-\a_k)(\ph(y_0)-\ph(y))]\,d\g_k
+\int_{\bG}(\ph(f(y,u))-\ph(y))\,d\g_k\\&\to \int_{\bG}(\ph(f(y,u))-\ph(y))\,d\g
=\int_{G}(\ph(f(y,u))-\ph(y))\,d\g \hbox{ as }k\to \infty.
\end{aligned}
\end{equation*}
Thus,
$$
\int_{G}(\ph(f(y,u))-\ph(y))\,d\g=0,
$$
that is, $\g\in W$. The proposition is proved. \hf

\bigskip

Denote
\begin{equation}\label{CV5}
\bar\G_{\a}:=\bigcup_{y_0\in \bY}(\G_{\a}(y_0)) \hbox{ and }\bar\G(S):=\bigcup_{y_0\in \bY}(\G(S,y_0)).
\end{equation}
Due to \eqref{f-inf}, any trajectory with initial condition at $\infty_Y$ stays at $\infty_Y$, and $(\infty_Y,u)$ can be identified with $\infty_G$ for any $u$. Therefore, $\G_{\a}(\infty_Y)=\{\delta_{\infty_G}\}$ (Dirac measure concentrated at $\infty_G$), and 
$$
\bar\G_{\a}=\G_{\a}\cup \{\delta_{\infty_G}\}.
$$
Similarly, $\G(S,y_{\infty})=\{\delta_{\infty_G}\}$ and $\bar\G(S)=\G(S)\cup \{\delta_{\infty_G}\}$ for all $S$. 

The following theorem follows from \cite{GPS-2017}, Theorem 5.4.
\begin{Theorem}
We have 
\begin{equation}\label{CV2}
\lim_{\a\to 1^-}(\bar{\rm co}\,\bar\G_{\a})=\lim_{S\to \infty}(\bar{\rm co}\,\bar\G({S,y_0}))=\bW.
\end{equation}
\end{Theorem}
{\bf Proof.} To apply Theorem 5.4 from \cite{GPS-2017} (that is, formula \eqref{CV1}) in the framework of compactified $Y$ and $G$ we need to make  the following changes:

(a) In \eqref{M17}, replace $Y$ with $\bY$ and $G$ with $\bG$. This leads to $W$ given by \eqref{M17} becoming $\bW$ given by \eqref{CV4}.

(b) In \eqref{CV1}, replace $\G_{\a}$ with $\bar\G_{\a}$ and $\G(S)$ with $\bar\G(S)$.

Then \eqref{CV1} becomes \eqref{CV2}. The theorem is proved. \hf

\bigskip

The following theorem ``decompactifies" relation \eqref{CV2} and provides a counterpart of \eqref{CV1} in the case of unbounded $Y$. 
\begin{Theorem}\label{Th4.3}
We have

{\rm (a)} \begin{equation}\label{M221}
(\liminf_{\a\to 1^-}({\rm co}\,\G_{\a}))\cap \P(G)=(\limsup_{\a\to 1^-}({\rm co}\,\G_{\a}))\cap \P(G)=W,
\end{equation}
{\rm (b)} \begin{equation}\label{M222}
(\liminf_{S\to \infty}({\rm co}\,\G(S)))\cap \P(G)=(\limsup_{S\to \infty}({\rm co}\,\G(S)))\cap \P(G)=W.
\end{equation}
\end{Theorem}

{\bf Proof.} We will only prove part (a), the proof of part (b) is similar.
Since  $\disp \limsup_{\a\to 1^-}({\rm co}\,\G_{\a})\subset \lim_{\a\to 1^-}(\bar{\rm co}\,\bar\G_{\a})$, from \eqref{CV2} and  \eqref{W-W} we conclude that 
\begin{equation}\label{M50}
(\limsup_{\a\to 1^-}({\rm co}\,\G_{\a}))\cap \P(G)\subset W.
\end{equation}

Let us show that $\disp W\subset\liminf_{\a\to 1^-}({\rm co}\,\G_{\a})$. Since $W\subset \P(G)$ due to \eqref{CV31}, together with \eqref{M50} this will imply \eqref{M221}.
Take $\g \in W$. From \eqref{CV2} and \eqref{W-W}
it follows that  for any $\a_k\to 1^-$ there exists a sequence  $\g_k\in \bar{\rm co}\,\bar\G_{\a_k}$ such that $\g_k\to \g$. Further, for each $k$ there exists a sequence $\g_{jk}\to \g_k$, $\g_{jk}\in {\rm co}\,\bar\G_{\a_k}$. By using a diagonalization argument, we can find a sequence $\g_{k}^*\in {\rm co}\,\bar\G_{\a_k}$ such that $\g_{k}^*\to \g$.

Since $ \bar\G_{\a}=\G_{\a}\cup\{\delta_{\infty_G}\}$, there exists $m=m(k)$ such that 
\begin{equation}\label{M30}
\g_k^*=\sum_{i=1}^{m} \l_{ik}\g_{ik}+\l_{m+1,k}\d_{\infty_G},\;\l_{ik}>0,\,\sum_{i=1}^{m+1} \l_{ik}=1,\,
\g_{ik}\in \G_{\a_k}.
\end{equation}
We can see that $\l_{m+1,k}\to 0$ as $k\to \infty$, since otherwise, if $\l_{m+1,k}\ge \beta>0$ along a subsequence,  for any $r>0$
$$
\int_{\bG} 1_{\bG\setminus rB}d\g^*_k\ge \beta,
$$
where  $B$ is the open unit ball in $\reals^m\times U_0$.
Hence,
$$
\int_{\bG} 1_{\bG\setminus rB}d\g\ge \beta,
$$
and
$$
\g(\{\infty_G\})=\lim_{r\to \infty}\int_{\bG} 1_{\bG\setminus rB}d\g\ge \beta
$$
contradicting the assumption that $\g$  is supported on $G$ ($\g\in W$). 

Set $\tilde \g_k$ to be the normalized first summation in \eqref{M30}, that is, 
$$
\tilde \g_k:={\sum_{i=1}^{m(k)} \l_{ik}\g_{ik}\o \sum_{i=1}^{m(k)}\l_{ik}}\in  {\rm co}\,\G_{\a_k}.
$$
Let us show that $\tilde \g_k\to \g$. Indeed, for any $q\in C(\bG)$ we have
\begin{equation}\label{CV10}
\int_{\bG} q\,d\tilde\g_k={\sum_{i=1}^{m(k)} \l_{ik}\int_{\bG}q\,d\g_{ik}\o \sum_{i=1}^{m(k)}\l_{ik}}
={1\o \sum_{i=1}^{m(k)}\l_{ik}}\left(\int_{\bG}q\,d\g_{k}^*-\l_{m+1,k}q(\infty_G)\right).
\end{equation}
Taking into account that, as $k\to \infty$,
$$
\l_{m+1,k}\to 0,\;\sum_{i=1}^{m(k)}\l_{ik}\to 1,\;\int_{\bG}q\,d\g_{k}^*\to\int_{\bG}q\,d\g,
$$
we conclude from \eqref{CV10} that $\disp \int_{\bG} q\,d\tilde\g_k\to \int_{\bG} q\,d\g$. Thus, for arbitrary $\a_k\to 1^-$ we constructed a sequence 
$\tilde \g_k\to \g$, $\tilde \g_k\in  {\rm co}\,\G_{\a_k}$, therefore, $\disp\g\in \liminf_{\a\to 1^-}({\rm co}\,\G_{\a})$.
The theorem is proved. \hf

\bigskip

{\bf Example.} Consider the one-dimensional system 
\begin{equation*}
\begin{aligned}
&y(t+1)=y(t)+u(t), \; t=0,1,\dots\, \\
&y(0)=y_0,\\
&y(t)\in Y=[0,\infty),\\
&u(t)\in U=\{0,1\}.
\end{aligned}
\end{equation*}
Admissible trajectories either run to infinity or reach a certain point and remain there. We have
 \begin{equation*}\label{M223}
W:=\{\g\in \P(Y\times U)|\, \int_{Y\times U}(\ph(y+u)-\ph(y))\,
\gm=0\quad \forall \ph\in C(\bY)\}.
 \end{equation*}
Let us verify that in this example $W=\P(Y\times \{0\})$, that is, $W$ consists of all probability measures supported on $Y\times \{0\}$. Indeed, it can be readily verified that for any $\g\in \P(Y\times \{0\})$ and $\ph\in C(\bY)$ we have $\disp\int_{Y\times U}(\ph(y+u)-\ph(y))\,\gm=0$. On the other hand, if $\g\in \P(Y\times U)$ is such that  $\g(Y\times \{1\})>0$, such measure doesn't belong to $W$, since for a monotonically increasing $\ph\in C(\bY)$ we have $\disp\int_{Y\times U}(\ph(y+u)-\ph(y))\,\gm>0$. 

Due to Theorem \ref{Th4.3}, we must have
$$
(\liminf_{\a\to 1^-}({\rm co}\,\G_{\a}))\cap \P(Y\times U)=(\limsup_{\a\to 1^-}({\rm co}\,\G_{\a}))\cap \P(Y\times U)=\P(Y\times \{0\})
$$
and
$$
(\liminf_{S\to \infty}({\rm co}\,\G(S)))\cap \P(Y\times U)=(\limsup_{S\to \infty}({\rm co}\,\G(S)))\cap \P(Y\times U)=\P(Y\times \{0\}).
$$
It can be intuitively understood why, for example, the inclusion 
\begin{equation}\label{ex1}
(\limsup_{S\to \infty}({\rm co}\,\G(S)))\cap \P(Y\times U)\subset \P(Y\times \{0\})
\end{equation}
must hold. If $\disp \g\in \limsup_{S\to \infty}({\rm co}\,\G(S))$ is such that $\g\notin \P(Y\times \{0\})$, that is, $\g(Y\times \{1\})>0$, then, as we will show, $\g\notin  \P(Y\times U)$. Since $\disp \g\in \limsup_{S\to \infty}({\rm co}\,\G(S))$, there exist sequences of times $S_k\to \infty$ and admissible processes $(y_k(t),u_k(t)),\,t=0,\dots,S_k-1$ that  generate occupational measures $\g_k$ such that $\g_k\to \g$. Since  $\g(Y\times \{1\})>0$ and $\g_k\to \g$, $u_k(t)=1$ occurs ``sufficiently frequently" to ensure that for any $y$ one has $\g_k([y,\infty)\times U)\ge\beta>0$ for sufficiently large $k$. Therefore, $\g([y,\infty)\times U)\ge\beta$, hence, $\g(\{\infty\}\times U)>0$, that is, $\g\notin  \P(Y\times U)$.

At the same time, the inclusion opposite to \eqref{ex1}, namely, $\disp\P(Y\times \{0\})\subset (\limsup_{S\to \infty}({\rm co}\,\G(S)))\cap \P(Y\times U)$, asserting that any measure in $\disp\P(Y\times \{0\})$ is a limit of measures from $({\rm co}\,\G(S))$, follows from Theorem \ref{Th4.3}, but is not obvious.

\section{Generalization of \eqref{eq-Res2}}\label{S-b}

In this section, we establish generalizations of \eqref{eq-Res2} to the situation when $Y$ is unbounded.

\subsection{The Property of Weakly$^*$ Convergent Sequences}

If $\g_k$ converges to $\g$ weakly$^*$ on $\P(\bG)$, then for any $q\in C(\bG)$ we have $\disp \lim_{k\to \infty}\int_{\bG}q(y,u)\,d\g_k\to \int_{\bG}q(y,u)\,d\g$ by definition of weak$^*$ convergence. As shown in the proposition below, this is also true when integration is taken over $G$ rather than $\bG$ as long as $\g_k$ and $\g$ are supported on $G$, and $q$ is a continuous bounded function, possibly without a limit at infinity. 

\begin{Proposition}\label{P-convergence}
Let $\g_k\in \P(G)$ and $\g\in \P(G)$ be such that $\g_k\to \g$ weakly$^*$ on $\P(\bG)$.  Then for any continuous bounded $q:\,G\to \reals$  we have
\begin{equation}\label{M36}
\lim_{k\to \infty}\int_{G}q(y,u)\,d\g_k= \int_{G}q(y,u)\,d\g.
\end{equation}
\end{Proposition}

{\bf Proof.} 
Let us see first that for any $\ve>0$ there exists $r'$ such that for all $r>r'$ we have $\disp \int_{G\setminus rB}1\,d\g_k\le \ve$ for all $k$. 

Indeed, assume it's not true. Then there exists $\ve_0>0$ such that for any $r>0$ there exists a subsequence (we do not relabel) such that $\disp \int_{G\setminus rB}1\,d\g_{k}> \ve_0$. Take $r$ such that $\disp \int_{G\setminus rB}1\,d\g< {\ve_0\o 2}$; such $r$ exists because the contrary would mean that $\g(\{\infty_G\})>0$ due to the property of measure continuity, while $\g$ is supported on $G$. 
Then
$$
\ve_0< \limsup_{k\to \infty}\int_{G\setminus rB}1\,d\g_{k}=\limsup_{k\to \infty}\int_{\bG\setminus rB}1\,d\g_{k}\le \int_{\bG\setminus rB}1\,d\g=\int_{G\setminus rB}1\,d\g<{\ve_0\o 2},
$$
which is a contradiction. 

Take $\ve>0$ and select $r$ so that 
\begin{equation}\label{M31}
\int_{G\setminus rB}1\,d\g\le \ve \hbox{ and }\int_{G\setminus rB}1\,d\g_k\le \ve\hbox{ for all }k.
\end{equation}
Due to the Tietze Extension Theorem, a continuous real-valued function can be extended from a closed subset of a normal topological space to the whole space without increasing its sup-norm. Denote by $\bar B$ the closed unit ball in $\reals^m\times U_0$ and apply this theorem to extend the function
\begin{equation*}
\begin{cases} q(y,u), &(y,u)\in G\cap r\bar B,\\ 0, &(y,u)\in G\setminus(r+1)B\end{cases}
\end{equation*}  
to $\reals^m\times U_0$ so that the extension, denoted  $q^*(y,u)$, satisfies 
$$
\disp \sup_{(y,u)\in \reals^m\times U_0}|q^*(y,u)|=\sup_{(y,u)\in G\cap r\bar B}|q(y,u)|\le\disp \sup_{(y,u)\in G}|q(y,u)|=:M.
$$ 
The function $q^*$ belongs to $C(\bG)$ by construction and we have
\begin{equation*}
\begin{aligned}
\left|\int_{G} q\,d\g_k-\int_{G} q\,d\g\right|\le \left|\int_{G} q\,d\g_k-\int_{G} q^*\,d\g_k\right|+\left|\int_{G} q^*\,d\g_k-\int_{G} q^*\,d\g\right|+
 \left|\int_{G} q^*\,d\g-\int_{G} q\,d\g\right|.
\end{aligned}
\end{equation*}
For the first difference we have   
$$
 \left|\int_{G} q\,d\g_k-\int_{G} q^*\,d\g_k\right|=\left|\int_{G\setminus rB} (q-q^*)\,d\g_k\right|\le  \int_{G\setminus rB} 2M\,d\g_k \le 2M\ve
$$
due to \eqref{M31}. The same estimate holds for the third difference. 
The second difference is equal to $ \left|\int_{\bG} q^*\,d\g_k-\int_{\bG} q^*\,d\g\right|$ and can be made arbitrarily small by increasing $k$. Therefore,\\ $ \left|\int_{G} q\,d\g_k-\int_{G} q\,d\g\right|\to 0$, that is, \eqref{M36} holds. The proposition is proved. \hf

\subsection{The Case when the Cost Function $g$ is Bounded}\label{S-bounded}

As in \eqref{M22}, let
\begin{equation*}
g^*:=\inf_{\g\in W} \int_{G} g(y,u)\,d\g.
\end{equation*}
Due to Proposition \ref{P-nonempty}, $g^*<\infty$ if there exists at least one trajectory of \eqref{A1} that stays in a bounded set for all $t$.

In this subsection, we consider the case when $g$ is {\em bounded} and derive generalizations of \eqref{eq-Res2} to the case of unbounded $Y$.

\begin{Proposition}\label{P-lower-bound-g} 
If $g$ is  bounded then 

{\rm (a)} \begin{equation}\label{M351}
g^*\ge \limsup_{\a\to 1^-}\inf_{y\in Y}V_{\a}(y),
\end{equation}
{\rm (b)} \begin{equation}\label{M352}
g^*\ge \limsup_{S\to \infty}\inf_{y\in Y}V(S,y).
\end{equation}
\end{Proposition}
{\bf Proof.} We prove part (a), the proof of part (b) is similar. Due to \eqref{CV11} we have
\begin{equation}\label{M37}
 \limsup_{\a\to 1^-}\inf_{y\in Y}V_{\a}(y)=\limsup_{\a\to 1^-}\inf_{\g\in \G_{\a}}\int_G g\,d\g=\limsup_{\a\to 1^-}\inf_{\g\in \rm{co}\,\G_{\a}}\int_G g\,d\g.
\end{equation}
From \eqref{M221} it follows that
\begin{equation}\label{M33}
g^*= \inf_{\g} \int_G g\,d\g,
\end{equation}
where inf is taken over  $\disp\g\in (\liminf_{\a\to 1^-}({\rm co}\,\G_{\a}))\cap \P(G)$.

Take $\d>0$ and $\disp\g_{\d}\in (\liminf_{\a\to 1^-}({\rm co}\,\G_{\a}))\cap \P(G)$ such that
\begin{equation}\label{M34}
\int_G g\,d\g_{\d}\le \inf_{\g}\int_G g\,d\g+\d,
\end{equation}
where, as before, inf  is taken over  $\disp\g\in (\liminf_{\a\to 1^-}({\rm co}\,\G_{\a}))\cap \P(G)$.
 Take any sequence $\a_k\to 1^-$. There exists a sequence $\g_k\to \g_{\d}$, $\g_k\in   {\rm co}\,\G_{\a_k}$  for which we have
\begin{equation*}
\int_G g\,d\g_k\to \int_G g\,d\g_{\d} \hbox{ as } k\to \infty
\end{equation*}
due to Proposition \ref{P-convergence}.  Therefore, for sufficiently large $k$ we have from \eqref{M34}
$$
\int_G g\,d\g_k\le \inf_{\g}\int_G g\,d\g+2\d,
$$
hence,
$$
\inf_{\g\in {\rm co}\,\G_{\a_k}} \int_G g\,d\g  \le \inf_{\g}\int_G g\,d\g+2\d,
$$
and
$$
\limsup_{\a\to 1^-} \inf_{\g\in {\rm co}\,\G_{\a}}\int_G g\,d\g\le \inf_{\g}\int_G g\,d\g=g^*.
$$
Taking into account this relation and \eqref{M37} we conclude that
$$
g^*\ge \limsup_{\a\to 1^-} \inf_{\g\in {\rm co}\,\G_{\a}}\int_G g\,d\g=\limsup_{\a\to 1^-}\inf_{y\in Y}V_{\a}(y).
$$
The proposition is proved. \hf

\bigskip

To obtain the inequalities opposite to \eqref{M351} and \eqref{M352}, introduce the following assumption.

(A1) {\em Optimal processes in the problems
\begin{equation}\label{A112}
\min_{y_0,\,u(\cdot)\in \U(y_0)}\sum_{t=0}^{\infty} \a^t g(y(t),u(t)) \hbox{\rm\; and }\min_{y_0,\,u(\cdot)\in \U_S(y_0)}\sum_{t=0}^{S-1} g(y(t),u(t)),
\end{equation} 
where minimization is taken with respect to both control and the initial condition, exist and remain in a bounded set $D\subset G$ for all $\a$ and $S$.}

\begin{Proposition}\label{P-upper-bound-g} 
Assume that $g$ is bounded and (A1) holds. Then 

{\rm (a)} 
\begin{equation}\label{M35}
g^*\le \liminf_{\a\to 1^-}\min_{y\in Y}V_{\a}(y),
\end{equation}
{\rm (b)} 
\begin{equation*}
g^*\le \liminf_{S\to\infty}\min_{y\in Y}V(S,y).
\end{equation*}

\end{Proposition}

{\bf Proof.} We prove part (a), the proof of part (b) is similar. 
Take an arbitrary sequence $\a_k\to 1^-$ and let $\g_k\in \G_{\a_k}$ be the occupational measures generated by optimal processes in \eqref{A112} with $\a=\a_k$. Then
\begin{equation}\label{M42}
\int_G g\,d\g_k=\min_{\g\in \G_{\a_k}}\int_G g\,d\g=\min_{\g\in {\rm co}\,\G_{\a_k}}\int_G g\,d\g.
\end{equation}
Take a subsequence of $\{\a_k\}$ along which $\disp\liminf_{\a_k\to 1^-}\min_{\g\in {\rm co}\,\G_{\a_k}}\int_G g\,d\g$  is reached. (We do not relabel.) Along a further subsequence, the measures $\g_k$ converge to some $\g^*\in \P(\bG)$ and, since $\g_k$ are supported on $D$, so is $\g^*$. (The proof of the latter fact is provided at the beginning of the proof of Proposition \ref{P-nonempty}.)
Due to Proposition \ref{P-convergence}, $\int_G g\,d\g_k\to \int_G g\,d\g^*$.   Therefore, from \eqref{M42} we have
\begin{equation}\label{M421}
\int_G g\,d\g^*= \liminf_{\a\to 1^-}\min_{\g\in {\rm co}\,\G_{\a}}\int_G g\,d\g.
\end{equation}
From \eqref{M221} it follows that
$$
g^*=\inf_{\g}\int_G g\,d\g,
$$
where inf on the right side is taken over $\disp\g \in (\limsup_{\a\to 1^-}({\rm co}\,\G_{\a}))\cap \P(G)$. Since $\disp\g^*\in (\limsup_{\a\to 1^-}({\rm co}\,\G_{\a}))\cap \P(G)$ and due to \eqref{M421}, we have 
\begin{equation}\label{M422}
g^*\le \int_G g\,d\g^*= \liminf_{\a\to 1^-}\min_{\g\in {\rm co}\,\G_{\a}}\int_G g\,d\g=\liminf_{\a\to 1^-}\min_{y\in Y}V_{\a}(y).
\end{equation}
The proposition is proved. \hf

\bigskip

As an immediate corollary of Propositions \ref{P-lower-bound-g} and \ref{P-upper-bound-g}, we obtain the following conditions ensuring \eqref{eq-Res2} in the case of unbounded $Y$:
\begin{Theorem}\label{Th} If $g$ is bounded and (A1) holds, then the limits $\disp\lim_{\a\to 1^-}\min_{y\in Y}V_{\a}(y)$ and $\disp\lim_{S\to \infty}\min_{y\in Y}V(S,y)$  exist and 
\begin{equation}\label{M38}
g^*= \lim_{\a\to 1^-}\min_{y\in Y}V_{\a}(y)=\lim_{S\to \infty}\min_{y\in Y}V(S,y).
\end{equation}
\end{Theorem}

Also note that we have proved the following:
\begin{Proposition}\label{Rem}
If $g$ is bounded and (A1) holds, then there exists a measure $\g^*$ supported on $D$ such that
$\disp g^*=\int_G g\,d\g^*$.
\end{Proposition}
{\bf Proof.} Since the inequality on the left side of \eqref{M422} holds as equality due to \eqref{M38}, the validity of the proposition follows. \hf

\bigskip

{\bf Example.} Take the system that we considered at the end of Section \ref{S-a}:
\begin{equation*}
\begin{aligned}
&y(t+1)=y(t)+u(t), \; t=0,1,\dots\, \\
&y(0)=y_0,\\
&y(t)\in Y=[0,\infty),\\
&u(t)\in U=\{0,1\}.
\end{aligned}
\end{equation*}
Let $g=g(y)$ in \eqref{A111} and \eqref{A112-1} be a bounded function of one variable $y$ that has a strict minimum at $\tilde y\ge 0$. It is clear that (A1) holds, the optimal process has initial condition $y_0=\tilde y$ and zero control, and $\disp\min_{y\in Y}V_{\a}(y)=\min_{y\in Y}V(S,y)=g(\tilde y)$ for all $\a$ and $S$. It is also clear that the minimum in $\disp \min_{\g\in W}\int_{Y\times U}g(y)\,d\g$ is reached at $\g$ being the Dirac function concentrated at the point where $y=\tilde y$ and $u=0$, which implies that   $\disp g^*:=
\min_{\g\in W}\int_{Y\times U}g(y)\,d\g=g(\tilde y)$. Thus,   $g^*=\disp\min_{y\in Y}V_{\a}(y)=\min_{y\in Y}V(S,y)$, as asserted in \eqref{M38}.

\subsection{The Case when the Cost Function $g$ is Unbounded}

If $g$ is unbounded from above (boundedness from below is assumed throughout), we truncate it to reduce to the framework of Section \ref{S-bounded}.

For $M>0$ denote
$$
g^M(y,u)=\min\{g(y,u),M\}.
$$
We introduce the following assumption, which is a strengthened version of (A1):

(A2) {\em There exists $M_1$ such that for all $M\ge M_1$ optimal processes in the problems
\begin{equation*}
\min_{y_0,\,u(\cdot)\in \U(y_0)}\sum_{t=0}^{\infty} \a^t g^M(y(t),u(t))  \hbox{\rm\; and }\min_{y_0,\,u(\cdot)\in \U_S(y_0)}\sum_{t=0}^{S-1} g^M(y(t),u(t)),
\end{equation*} 
where minimization is taken with respect to both control and the initial condition, exist and remain in a bounded set $D\subset G$ for all $\a$ and $S$.} 

\begin{Lemma} 
If (A2) holds then so does (A1), and for all $\a,S$ and sufficiently large $M$ we have 
\begin{equation}\label{CV9}
\min_{y\in Y}V_{\a}^M(y)=\min_{y\in Y}V_{\a}(y),
\end{equation}
\begin{equation}\label{CV91}
\min_{y\in Y}V^M(S,y)=\min_{y\in Y}V(S,y),
\end{equation}
where $\disp V_{\a}^M(y_0)=\min_{u(\cdot)\in \U(y_0)}(1-\a)\sum_{t=0}^{\infty} \a^t g^M(y(t),u(t))$ and \\$\disp V^M(S,y_0)=\min_{u(\cdot)\in \U_S(y_0)}{1\o S}\sum_{t=0}^{S-1} \a^t g^M(y(t),u(t))$.
\end{Lemma}

{\bf Proof.} We will prove \eqref{CV9}, the proof of \eqref{CV91} is similar. For an admissible process $(y(\cdot),u(\cdot))$ denote
$$
J_{\a}(u(\cdot)):=(1-\a)\sum_{t=0}^{\infty} \a^t g(y(t),u(t)) \hbox { and }J_{\a}^M(u(\cdot)):=(1-\a)\sum_{t=0}^{\infty} \a^t g^M(y(t),u(t)).
$$
Assume that (A1) does not hold for the problem with discounting, that is, for some $\a$ there exists an admissible process $(\by(\cdot),\bu(\cdot))$, not contained in $D$, such that
$J_{\a}(\bu(\cdot))\le J_{\a}(u(\cdot))$ for any process $(y(\cdot),u(\cdot))$ contained in $D$. For any $M$ we have $J_{\a}^M(\bu(\cdot))\le J_{\a}(\bu(\cdot))$ and for $\disp M\ge \max_{(y,u)\in D}g(y,u)$ we have $J_{\a}(u(\cdot))= J_{\a}^M(u(\cdot))$. Putting these together, we get 
$$
J_{\a}^M(\bu(\cdot))\le J_{\a}(\bu(\cdot))\le J_{\a}(u(\cdot))= J_{\a}^M(u(\cdot)),
$$
which contradicts (A2) for the problem with discounting for $M\ge M_1$. Thus, (A1) holds with $\disp M_0=\max\{M_1,\max_{(y,u)\in D}g(y,u)\}$. Since $g^M=g$ on $D$ for $M\ge M_0$, \eqref{CV9} follows.
The lemma is proved. \hf

\bigskip

\begin{Theorem}\label{Th3}  Assume (A2). Then the limits $\disp \lim_{\a\to 1^-}\min_{y\in Y}V_{\a}(y)$ and $\disp\lim_{S\to \infty}\min_{y\in Y}V(S,y)$ exist and
\begin{equation*}
g^*= \lim_{\a\to 1^-}\min_{y\in Y}V_{\a}(y)=\lim_{S\to \infty}\min_{y\in Y}V(S,y).
\end{equation*}
\end{Theorem}

{\bf Proof. }  We will prove the first equality.
Since (A2) implies (A1) with $\disp M_0=\max\{M_1,\max_{(y,u)\in D}g(y,u)\}$ and $g^M$ is bounded, from Theorem \ref{Th} we have for any $M\ge M_0$
\begin{equation}\label{M23}
\lim_{\a\to 1^-} \min_{y\in Y}V_{\a}^M(y)=(g^M)^*,
\end{equation}
where 
\begin{equation}\label{M24}
(g^M)^*:=\inf_{\g\in W} \int_{G} g^M(y,u)d\g.
\end{equation}
Due to Proposition \ref{Rem}, the measure $\g^*\in W$ such that 
$\disp (g^M)^*=\int_{G} g^M(y,u)d\g^*$ exists and is supported on $D$. 
Since $g^M(y,u)=g(y,u)$ for all $(y,u)\in D$ if $\disp M\ge \max_{(y,u)\in D}g(y,u)$, we have
\begin{equation}\label{M27}
(g^M)^*=\int_{G} g^M(y,u)d\g^*=\int_{G} g(y,u)d\g^*\ge g^*.
\end{equation}
Since it is obvious that $(g^M)^*\le g^*$, we conclude that the last inequality holds as equality.  
The statement of the theorem now follows \eqref{M23} and \eqref{CV9}. The theorem is proved. \hf

\bigskip

{\bf Acknowledgement.} The author wishes to express his gratitude to V. Gaitsgory for helpful suggestions during this research.

Email address of the author: ius13@psu.edu


\begin{thebibliography}{99}


\bibitem{Ash} R. Ash,
{Measure, Integration and Functional Analysis},
Academic Press, 1972.

\bibitem{BG05} V. Borkar, V. Gaitsgory, On Existence of Limit Occupational Measures Set of a Controlled Stochastic Differential Equation, SIAM J. Control Optim., 44(4) (2005), 1436-1473.

\bibitem{BG}
V. Borkar, V. Gaitsgory, Linear Programming Formulation of Long
Run Average Optimal Control Problem,  {  J. of Optimization Theory and Applications},   181(1)  (2019), 101--125.

\bibitem{BGS}
V. Borkar, V. Gaitsgory and I. Shvartsman,
 LP Formulations of  Discrete Time Long-Run Average Optimal Control Problems: The Non-Ergodic Case , {  SIAM Journal on Control and Optimization},   57(3)  (2019), 1783--1817.



\bibitem{BQR-2015} R.\ Buckdahn, M.\ Quincampoix and J.\ Renault,   On Representation Formulas for Long Run Averaging Optimal Control Problem,  {  Journal of Differential Equations},   259(11)  (2015), 5554--5581.



\bibitem{G04} V.\ Gaitsgory, On a Representation of the Limit Occupational Measure of a Control System with Applications to Singularly Perturbed Control Systems, {  SIAM J.\ of Control and Optimization},    43(1) (2004), 325--340.


\bibitem{GQ} V.\ Gaitsgory and M.\ Quincampoix,
 Linear programming approach to deterministic
infinite horizon optimal control problems with discounting,   {  SIAM J.\ of Control and
Optimization},    48(4) (2009), 2480--2512.



\bibitem {GPS-2017} V.\ Gaitsgory, A.\ Parkinson and I.\ Shvartsman,  Linear programming formulations of deterministic infinite horizon optimal control problems in discrete time,  {  Discrete and Continuous Dynamical Systems Series B},    22(10) (2017),  3821--3338.


\bibitem{GPS-2019}  V.\ Gaitsgory, A.\ Parkinson and I.\ Shvartsman,   Linear programming based optimality conditions and approximate solution of a deterministic infinite horizon discounted optimal control problem in discrete time,  {  Discrete and Continuous Dynamical Systems, Series B},   24(4)  (2019), 1743--67.

\bibitem{GS-Opt-22}
V. Gaitsgory, I. Shvartsman, LP-Related Representations of Cesaro and Abel Limits of Optimal Value Functions, Optimization, 71(4) (2022), 1151-1170.

\bibitem{GS-DCDS-22}
V. Gaitsgory, I. Shvartsman, Linear Programming Estimates for Cesaro and Abel Limits of Optimal Values in Optimal Control Problems, Discrete and Continuous Dynamical Systems, Series B, 27(3) (2022), 1591--1610.






\bibitem{HLerma23} O. Hern\'andez-Lerma, L. Laura-Guarachi, S. Mendoza-Palacios, A survey of average cost problems in deterministic discrete-time control systems, J. of Mathematical Analysis and Applications, 522(1) (2023).




\bibitem{Sorin92}
 E.\ Lehrer and S.\  Sorin,
   A uniform Tauberian theorem in dynamic programming, 
{  Mathematics of Operations Research},   17(2)  (1992), 303--307.

\bibitem{OV-2012}  M.\ Oliu-Barton and G.\   Vigeral,   A uniform Tauberian theorem in optimal control, in {  ``Annals of International Society of Dynamic Games"} (eds. P.\ Cardaliaguet and R.\  Grossman),   12,   199--215, Birkhauser/Springer, New York (2013).



\bibitem{ShvPAFA21} I. Shvartsman, Lack of Equality between Abel and Cesaro Limits in Discrete Optimal Control and the Implied Duality Gap, Pure and Applied Functional Analysis, 6(6) (2021), 1495--1507. 


\bibitem{ShvCOT} I. Shvartsman, Optimality Conditions in Discrete-Time Infinite-Horizon Optimal Control Problem with Discounting, Communications in Optimization Theory, Vol. 2023, 2023, pp. 1-10.






\end{thebibliography}
\end{document}